\magnification=1200 
\vsize=187mm 
\hsize=125mm 
\hoffset=4mm
\voffset=10mm

%XXXXXXXXXXXXXXXXXXXXXXXXXXXXXXXXXXXXXXXXXXXXXXXXXXXXXXXXXXXXXXXXXXXXXXXXXXXXXXXXXXXX

\abovedisplayskip=4.5pt plus 1pt minus 3pt
\abovedisplayshortskip=0pt plus 1pt
\belowdisplayskip=4.5pt plus 1pt minus 3pt
\belowdisplayshortskip=2.5pt plus 1pt minus 1.5pt
\smallskipamount=2pt plus 1pt minus 1pt
\medskipamount=4pt plus 2pt minus 1pt
\bigskipamount=9pt plus 3pt minus 3pt

%XXXXXXXXXXXXXXXXXXXXXXXXXXXXXXXXXXXXXXXXXXXXXXXXXXXXXXXXXXXXXXXXXXXXXXXXXXXXXXXXXXXXXXXX

%XXXXXXXXXXXXXXXXXXXXXXXXXXXXXXXXXXXXXXXXXXXXXXXXXXXXXXXXXXXXXXXXXXXXXXXXXXXXXXXXXXXXXXXX
%                   Les Hauts de Pages
%XXXXXXXXXXXXXXXXXXXXXXXXXXXXXXXXXXXXXXXXXXXXXXXXXXXXXXXXXXXXXXXXXXXXXXXXXXXXXXXXXXXXXXXX

\newif\ifpagetitre          \pagetitretrue
\newtoks\hautpagetitre     \hautpagetitre={\hfil}
\newtoks\baspagetitre     \baspagetitre={\hfil\tenrm\folio\hfil}

\newtoks\auteurcourant     \auteurcourant={\hfil}
\newtoks\titrecourant     \titrecourant={\hfil}
\newtoks\chapcourant     \chapcourant={\hfil}

\newtoks\hautpagegauche     \newtoks\hautpagedroite
\hautpagegauche={\vbox{\it\noindent\the\chapcourant\hfill\the\auteurcourant\hfill
{ }\smallskip\smallskip\vskip 2mm\line{}}}
\hautpagedroite={\vbox{\hfill\it\the\titrecourant\hfill{ }
\smallskip\smallskip\vskip 2mm\line{}}}

\newtoks\baspagegauche     \newtoks\baspagedroite
\baspagegauche={\hfil\tenrm\folio\hfil}
\baspagedroite={\hfil\tenrm\folio\hfil}

\headline={\ifpagetitre\the\hautpagetitre
\else\ifodd\pageno\the\hautpagedroite
\else\the\hautpagegauche\fi\fi}

\footline={\ifpagetitre\the\baspagetitre
\global\pagetitrefalse
\else\ifodd\pageno\the\baspagedroite
\else\the\baspagegauche\fi\fi}

%XXXXXXXXXXXXXXXXXXXXXXXXXXXXXXXXXXXXXXXXXXXXXXXXXXXXXXXXXXXXXXXXXXXXXXXXXXXXXXXXXXXXXXX

\auteurcourant={}
\titrecourant={}

%XXXXXXXXXXXXXXXXXXXXXXXXXXXXXXXXXXXXXXXXXXXXXXXXXXXXXXXXXXXXXXXXXXXXXXXXXXXXXXXXXXXXXX
%XXXXXXXXXXXXXXXXXXXXXXXXXXXXXXXXXXXXXXXXXXXXXXXXXXXXXXXXXXXXXXXXXXXXXXXXXXXXXXXXXXXXXX
%                                 MACROS PERSONNELLES
%XXXXXXXXXXXXXXXXXXXXXXXXXXXXXXXXXXXXXXXXXXXXXXXXXXXXXXXXXXXXXXXXXXXXXXXXXXXXXXXXXXXXXX
%XXXXXXXXXXXXXXXXXXXXXXXXXXXXXXXXXXXXXXXXXXXXXXXXXXXXXXXXXXXXXXXXXXXXXXXXXXXXXXXXXXXXX

\font\timeonze=cmr10 scaled 1100

\font\bfonze=cmbx10 scaled 1100

\def\NN{{\mathord{I\!\! N}}}
\def\RR{{\mathord{I\!\! R}}}
\def\EE{{\mathord{I\!\! E}}}

\def\CC{{\mathord{C\mkern-16mu{\phantom t\vrule}{\phantom o}}}}

\def\Rp{\RR^+}
\def\pro#1{{(#1_t)}_{t\geq 0}}

\def\normca#1{{\left\vert\left\vert #1\right\vert\right\vert}^2}
\def\ab#1{\left\vert #1\right\vert}

\def\titredeux#1#2#3{\centerline{\bfonze#1}\medskip
     \centerline{\bfonze#2}\bigskip\bigskip\bigskip\centerline{\timeonze#3}}

\def\spa#1#2{\bigskip\medskip\noindent{\bfonze #1\ #2}\par\nobreak\bigskip}

\def\th#1{\bigskip\smallskip\noindent{\bf Theorem #1}$\,$--$\,$}

\def\prp#1{\bigskip\smallskip\noindent{\bf Proposition #1}$\,$--$\,$}

\def\prf{\bigskip\noindent{\bf Proof}\par\nobreak\smallskip}

\def\findem{\hfill\hbox{\vrule height 2.5mm depth 0 mm width 2.5 mm}}
\def\indic{{\mathop{\rm 1\mkern-4mu l}}}

\def\eq{\eqalignno}
\def\hf{\hfill}
\def\wh{\widehat}
\def\wt{\widetilde}
\def\ps#1#2{{<}\, #1\, ,\, #2\,{>}}

\def\ld{\ldots}
\def\cd{\cdot}

\def\qed{\findem}
\def\thsp{\thinspace}

\def\cd{\cdot}

\def\aa{a.a.\thsp}

\def\frac#1#2{{{#1}\over{#2}}}

\def\rA{{\cal A}}\def\rB{{\cal B}}\def\rD{{\cal D}}
\def\rF{{\cal F}}
\def\rL{{\cal L}}
\def\rM{{\cal M}}\def\rP{{\cal P}}
\def\rS{{\cal S}}

\def\a{\alpha}

\def\d{\delta}
\def\e{\varepsilon}

\def\s{\sigma}

\def\o{\omega}

\def\G{\Gamma}
\def\D{\Delta}

\def\O{\Omega}

\def\F{\Phi}

\titredeux{FROM $(n+1)\,$-LEVEL ATOM CHAINS}{TO $n$-DIMENSIONAL
NOISES}{St\'ephane ATTAL and Yan PAUTRAT}
\bigskip\bigskip\bigskip
\centerline{\it This article is dedicated to the memory of
Paul-Andr\'e MEYER} 
\bigskip
\spa{}{Abstract}
{\sevenrm In quantum physics, the state space of a countable chain of
$(n+1)\,$-level atoms becomes, in the continuous field limit, a Fock space
with multiplicity $n$. In a  more functional analytic language, the
continuous tensor product space over $\Rp$ of copies of the space $\CC^{n+1}$ is the
symmetric Fock space $\G_s(L^2(\Rp;\CC^n))$. In this article we focus
on the probabilistic interpretations of these facts. We show that they
correspond to the approximation of the $n$-dimensional normal
martingales by means of obtuse random walks, that is, extremal random
walks in $\RR^n$ whose jumps take exactly $n+1$ different values. We
show that these probabilistic approximations are carried by the
convergence of the 
basic matrix basis $a^i_j(p)$ of $\otimes_\NN\CC^{n+1}$ to the usual
creation, annihilation and gauge processes on the Fock space.} 

\spa{I.}{Introduction}
In functional analysis, the tensor product of a family of
Hilbert spaces indexed by a continuous set, is a well-understood
notion (see the very complete book [Gui]) which leads to notions such as
``Fock spaces'' or  ``symmetric space associated to a measured space''. 

A  physical interpretation of  those continuous tensor
product spaces consists in considering them  as  the continuous
field limit of a countable chain of quantum system state spaces (such
as a spin chain, for example).

The interesting point in these constructions is that, for all
$n\in\NN$, the
continuous tensor product space
$$
\bigotimes _{\Rp}\CC^{n+1}
$$
is the symmetric Fock space $\G_s(L^2(\Rp;\CC^n))$. In a more physical
language, the continuous field limit of the state space of a countable
chain of 
$(n+1)\,$-level atoms is a Fock space with multiplicity $n$. A
rigourous setting in which such an approximation is made true is
developped in [At1].
\bigskip
Both the spaces
$\otimes_\NN\CC^{n+1}$ and 
$\G_s(L^2(\Rp;\CC^n))$ admit natural probabilistic
interpretations. Indeed, 
 the Fock space
$\G(L^2(\Rp;\CC^n))$ admits natural probabilistic interpretations in
terms of $n$-dimensional  normal martingales, such as $n$-dimensional
Brownian motion, $n$-dimensional Poisson process,  $n$-dimensional Az\'ema
martingales ... (cf [A-E] and [At2]). The aim of this article is to
understand 
how the approximation of $\G(L^2(\Rp;\CC^n))$ by means of spaces
$\otimes_\NN\CC^{n+1}$ can be interpreted in probabilistic terms. 

The structure of the space $\otimes_\NN\CC^{(n+1)}$
suggests that we are dealing with  random walks whose jumps are taking
$(n+1)$ different values. 

In this article we show that  the key point of this approximation is
the notion of {\it obtuse random walks}, developped in [A-E]. They are
the centered and normalized random variables in $\RR^n$ 
which take exactly $(n+1)$ different values. 

These obtuse random variables are naturally associated to an algebraic
object called {\it
sesqui-symmetric 3-tensor} and  the associated random walk
satisfies a {\it discrete-time structure equation}. This structure
equation allows us to represent the multiplication operators by this
random walk in terms of some basic operators of
$\otimes_\NN\CC^{n+1}$. 

Considering
the approximation of the Fock space  $\G(L^2(\Rp;\CC^n))$ by means of
spaces $\otimes_\NN\CC^{n+1}$, we obtain
the approximation of  a continuous-time
normal martingale. The sesqui-symmetric 3-tensor $\F$ then converges
to a so-called {\it doubly-symmetric 3-tensor} which is
the key of the structure equation describing the probabilistic
behaviour of that normal martingale (jumps, continuous and purely
discontinuous parts...). 
\bigskip
This article is organized in the following way: in section two we
introduce the state space of the atom chain and the associated
operators. In section three, we describe obtuse random walks in $\RR^n$, their
structure equations and their representations as operators on the state space of
the atom chains. In section four we introduce Fock space and its quantum
stochastic calculus, and the relation of these objects with the atom
chains. In section five we describe structure equations for normal
martingales and the information given by these equations in a special
case. In section six we put together all of our tools and prove
convergence in law of random walks to well-identified normal
martingales. In section seven we review some explicit and illustrative examples.

\spa{II.}{The structure of the atom chain}

We here introduce the mathematical structure and notations associated
to the space $\otimes_\NN\CC^{n+1}$. As the reader will easily see,
this only means choosing a particular basis for the vectors and  for
the operators on that space. The physical-like 
terminology that we here use time to time is not necessary for the
sequel, it is just informative (though it is pertinent and really
used in articles such as [A-P]).
\bigskip
Consider the space $\CC^{n+1}$ in which we choose an orthonormal basis
denoted by $\{\O,X^1,\ld,X^n\}$. This space and this particular choice
of an orthonormal basis physically represent either a particle with $n$
excited states $X^i$ and a ground state $\O$, or a site which is
either empty ($\O$) or occupied by a type $i$ particle ($X^i$). We often
write $X^0$ for $\O$ when we need unified notations, but it is
important in the sequel to distinguish one of the basis states. 

Together with this basis of $\CC^{n+1}$ we consider the following natural basis of $\rL(\CC^{n+1})=M_{n+1}(\CC)$:
$$
a^i_jX^k=\d_{ki}X^j,
$$
for all $i,j,k=0,\ld,n$. With these notations the operator $a^0_j$
corresponds, up to a sign factor, to classical fermionic creation operator for the particle
$X^j$; indeed, we have $a^0_j\O=X^j$ and $(a^0_j)^2=0$. The operator
$a^j_0$ corresponds to its associated annihilation operator. The
operator $a^i_j$ exchanges a $i$-level state with a $j$-level state
particle.
\smallskip
\def\TF{{\rm T}\!\F}
We now consider a  chain of copies of this system, like a chain of
$(n+1)$-level atoms. That is, we consider the
Hilbert space 
$$
\TF=\bigotimes_{i\in\NN}\CC^{n+1}
$$
 made of a countable tensor
product, indexed by $\NN$, of copies of $\CC^{n+1}$. By this we mean
 that a natural
orthonormal basis of $\TF$ is described by the family
$$
\left\{X_A;A\in\rP_n\right\}
$$
where

-- $\rP_n$ is the set of finite subset $A=\{(n_1,i_1),\ld,(n_k,i_k)\}$
   of
   $\NN\times\{1,\ld,n\}$ such that the $n_i$'s are two by two
   different. Another way to describe the set $\rP_n$ is to identify
   it to the set of sequences ${(A_k)}_{k\in\NN}$ with values in $\{0,\ld,n\}$,
   but taking only finitely many times a value different from 0.

-- \ $X_A$ denotes the vector
$$
\O\otimes\ld\otimes\O\otimes X^{i_1}\otimes
\O\otimes\ld\otimes\O\otimes X^{i_2}\otimes\ld
$$
of $\TF$, where $X^{i_1}$ appears in the copy number $n_1$, $X^{i_2}$
appears in the copy $n_2$,... When $A$ is seen as a sequence ${(A_k)}_{k\in\NN}$
as above, then $X_A$ is advantageously written $\otimes_{k}X_{A_k}$.
\smallskip
The physical meaning of this basis is easy to understand: we
have a chain of sites, indexed by $\NN$; on each site there is an atom
in the ground state or an atom  at energy level 1... The above basis
vector $X_A$ 
specifies that there is an atom at level $i_1$ in the site
$n_1$, an atom at level $i_2$ in the site $n_2\ld$, all the other
sites being at the ground state. The space $\TF$ is what we shall call
the {\it 
$(n+1)$-level atom
chain}.
\bigskip
We denote by $a^i_j(k)$ the natural ampliation of the operator $a^i_j$
to $\TF$ which acts as $a^i_j$ on the copy number $k$ of $\CC^{n+1}$
and as the identity 
on the other copies.  
\bigskip
Note, for information only,  that the operators $a^i_j(k)$ form a
basis of the algebra $\rB(\TF)$ of bounded operators on $\TF$. That
is, the von Neumann algebra 
generated by the $a^i_j(k)$, $i,j=0,\ld,n$, $k\in\NN$, is the whole
of $\rB(\TF)$ (for $\TF$  admits no subspace which is non trivial
and invariant under this algebra).

\spa{III.}{Obtuse random walks in $\RR^n$}
We now abandon for a while this structure in order to concentrate on
the probabilistic and algebraic structure of the obtuse random
variables. The space $\TF$ will come back naturally when
describing the obtuse random walks.
\bigskip
Let $X$ be a random variable in $\RR^n$ which takes exactly $n+1$
different values $v_1,\ld,v_{n+1}$ with respective probability
$\a_1,\ld,\a_{n+1}$ (all different from 0 by hypothesis). We assume,
for simplicity, that $X$ is defined on its canonical space
$(A,\rA,P)$, that is, $A=\{1,\ld,n+1\}$, $\rA$ is the $\s$-field of
subsets of $A$, the probability measure $P$ is  given by
$P(\{i\})=\a_i$ and $X$ is given by $X(\{i\})=v_i$, for all
$i=1,\ld,n+1$.
 
Such a random variable $X$ is called {\it centered and normalized} if
$\EE[X]=0$ and $\hbox{Cov}(X)=I$.

A family of elements $v_1,\ld,v_{n+1}$ of $\RR^n$ is called an {\it
obtuse system} if $$\ps{v_i}{v_j}=-1$$ for all $i\not =j$.

We consider the coordinates $X^1,\ld,X^n$ of $X$ in the canonical
basis of $\RR^n$, together with the random variable $\O$ on
$(A,\rA,P)$ which is deterministic and always equal to 1.

We put $\wt X^i$ to be the random variable $\wt
X^i(j)=\sqrt{\a_j}\, X^i(j)$ and $\wt\O(j)=\sqrt{\a_j}$. For any element $v=(a_1,\ld,a_n)$ of
$\RR^n$ we put $\wh v=(1,a_1,\ld,a_n)\in\RR^{n+1}$.

The following proposition is rather straightforward and left to the
reader.

\prp{1.}{\it The following assertions are equivalent.
\smallskip
i) $X$ is centered and normalized.
\smallskip
ii) The $(n+1)\times(n+1)$-matrix $(\wt\O,\wt X^1,\ld,\wt X^n)$ is
unitary.
\smallskip
iii)  The $(n+1)\times(n+1)$-matrix $(\sqrt{\a_1}\, \wh v_1,\ld,
\sqrt{a_{n+1}}\, \wh v_{n+1})$ is unitary.
\smallskip
iv) The family $v_1,\ld,v_{n+1}$ is an obtuse system of $\RR^n$ and 
$$
\a_i=\frac 1{1+\normca{v_i}}.
$$}

Let $T$ be a 3-tensor in $\RR^n$, that is, a linear mapping from
$\RR^n$ to $M_n(\RR)$. We write $T^{ij}_k$ for the coefficients of $T$ in
the canonical basis of $\RR^n$, that is,
$$
(T(x))_{i,j}=\sum_{k=1}^nT^{ij}_kx_k.
$$
Such a 3-tensor $T$ is called {\it sesqui-symmetric} if 
\smallskip
i) $(i,j,k)\longmapsto T^{ij}_k$ is symmetric

\noindent and

ii) $(i,j,l,m)\longmapsto\sum_kT^{ij}_kT^{lm}_k+\d_{ij}\d_{lm}$ is
symmetric.
\bigskip
\th{2}{\it If $X$ is a centered and normalized random variable in
$\RR^n$, taking exactly $n+1$ values, then there exists a
sesqui-symmetric 3-tensor $T$ such that 
$$
X\otimes X=I+T(X).\eqno{(1)}
$$
}
\prf
By Proposition 1, the matrix $(\sqrt{\a_1}\,\,\wh v_1,\ld,
\sqrt{\a_{n+1}}\,\, \wh v_{n+1})$ is unitary. In particular the matrix
$(\wh v_1,\ld, 
\wh v_{n+1})$ is invertible and so is its adjoint matrix. But the latter
is the matrix whose columns are the values of the random variables
$\O,X_1,\ld,X_n$. As a consequence, these $n+1$ random variables are
linearly independent. They thus form a basis of $L^2(A,\rA,P)$ for it
is a $n+1$ dimensional space. 

The random variable $X^iX^j$ belongs to  $L^2(A,\rA,P)$ and can thus be
written as
$$
X^iX^j=\sum_{k=0}^nT^{ij}_kX^k,
$$
for some real coefficients $T^{ij}_k$,
$k=0,\ld,n$, $i,j=1,\ld n$, where $X^0$ denotes $\O$. The fact that $\EE[X^k]=0$ and
$\EE[X^iX^j]=\d_{ij}$ implies $T^{ij}_0=\d_{ij}$. This gives the
representation (1).

The fact that the 3-tensor $T$ associated to the above coefficients
$T^{ij}_k$, $i,j,k=1,\ld n$, is sesqui-symmetric is an easy
consequence of the fact that the expressions $X^iX^j$ are symmetric in
$i,j$ and $X^i(X^jX^k)=(X^iX^j)X^k$ for all $i,j,k$. We leave it to
the reader.\qed
\bigskip
There is actually a natural bijection between the set of sesqui-symmetric  3-tensors
and the set of obtuse random variables. This is a result  obtained in [A-E],
Theorem 2, pp. 268-272, which is far from obvious but which we shall
not really need here. 

\th{3.}{\it The formulas
$$
S=\{x\in\RR^n; x\otimes x=I+T(x)\}.
$$
and 
$$
T(x)=\sum_{y\in S}p_y\ps{y}{x}\,y\otimes y,
$$
where $p_x=1/(1+\normca x)$, define a bijection between the set of
sesqui-symmetric 3-tensor $T$ on $\RR^n$ and the set of obtuse systems
$S$ in $\RR^n$.}\qed 
\bigskip
Now we wish to consider the random walks which are induced by obtuse
systems. That is, on the probability space 
$(A^\NN,\rA^{\otimes\NN},P^{\otimes \NN})$, we consider a sequence
${(X(n))}_{n\in\NN}$ of independent random variables with the
same law as a given 
centered normalized random variable $X$.

Recalling the notations of section II, for any $A\in\rP_n$, we define the random variable
$$
X_A=\prod_{(p,i)\in A}X^i(p)
$$
with the convention
$$
X_\emptyset=\indic.
$$
\def\esp#1{\EE[#1]}
\prp{4.}{\it The family $\{X_A;A\in\rP_n\}$ is an orthonormal basis
of the space
$L^2(A^\NN,\rA^{\otimes\NN},P^{\otimes \NN})$.}

\prf
For any $A,B\in\rP_n$ we have
$$
\ps{X_A}{X_B}=\esp{X_A X_B}=\esp{X_{A\D B}}\esp{X_{A\cap B}^2}
$$
by the independence of the $X(p)$. For the same reason, the first term
$\esp{X_{A\D B}}$ 
gives 0 unless $A\D B=\emptyset$, that is $A=B$. The second term
$\esp{X_{A\cap B}^2}$ is then equal to $\prod_{(p,i)\in
A}\esp{X^i(p)^2}=1$. This proves the orthonormal character of the
family  $\{X_A;A\in\rP_n\}$.

Let us now prove that it generates a dense subspace of
$L^2(A^\NN,\rA^{\otimes\NN},P^{\otimes \NN})$. Had we considered
random walks indexed by $\{0,\ld, N\}$ instead of $\NN$, the $X_A$, $A\subset\{0,\ld, N\}$ would have formed an
orthonormal basis of $L^2(A^N,\rA^{\otimes N},P^{\otimes N})$, for their
dimensions are equal. Now a general element $f$ of
$L^2(A^\NN,\rA^{\otimes\NN},P^{\otimes \NN})$ can be easily approached
by a sequence $(f_N)_N $ such that  $f_N\in L^2(A^N,\rA^{\otimes N},P^{\otimes N})$,
for all $N$, by taking conditional expectations on the trajectories
of 
$X$ up to time $N$.  \qed
\bigskip
For every obtuse random variable $X$, we thus obtain a Hilbert space 
$$\TF(X)=L^2(A^\NN,\rA^{\otimes\NN},P^{\otimes \NN}),$$ with a natural
orthonormal basis $\{X_A;A\in\rP_n\}$ which emphasizes the independence of
the $X(p)$'s. In particular there is a natural isomorphism between all
the spaces $\TF(X)$ which consists in identifying the associated
bases. In the same way, all these canonical spaces $\TF(X)$ of obtuse random
walks are naturally isomorphic to the atom chain $\TF$ of previous
section (again by identifying their natural orthonormal bases). 

Of course this identification of Hilbert spaces does not mean much for the
moment: in particular, it loses all the probabilistic properties of
the random variables $X^i(p)$, be it individual (the law) or
collective (probabilistic independence) properties.

The only way to recover the full probabilistic information on $X^i(p)$
in the Hilbert space formalism associated to $\TF$ is to consider the
{\it multiplication operator} by  $X^i(p)$ instead of the Hilbert
space element $X^i(p)$. Indeed, if we know the representation in $\TF$
of the operator $\rM_{X^i(p)}$ of multiplication  by $X^i(p)$ on
$\TF(X)$, we know everything about the random variable $X^i(p)$ and its
relation with other random variables. The above idea is what makes
quantum probabilistic tools relevent for the study of classical
probability; following this idea, the next theorem is one of the keys
of this article. It is what allows us to translate probabilistic
properties into theoretic language, showing that 
{\it all} the obtuse random walks in $\RR^n$
can be represented in a single space $\TF$ with very economical means:
linear combinations of the operators $a^i_j(p)$.

\th{5.}{\it Let $X$ be an obtuse random variable, let
${(X(p))}_{p\in\NN}$ be the associated random walk on the
canonical space $\TF(X)$. Let $T$ be the sesqui-symmetric 3-tensor
associated to $X$. Let $U$ be the natural unitary isomorphism from
$\TF(X)$ to $\TF$. Then, for all $p\in\NN,i=\{1,\ld,n\}$ we have
$$
U\rM_{X_i(p)}U^\ast=a^0_i(p)+a^i_0(p)+\sum_{j,l=1}^n T^{jl}_i\,a^j_l(p).
$$
}

\prf
It suffices to compute the action of $\rM_{X_i(p)}$ on the basis elements
$X_A$, $A\in\rP_n$. Denote by ``$(p,.)\not\in A$'' the claim ``for no
$i$ does $(p,i)$ belong to $A$''. Then, by Theorem 1, there exists a
sesquisymmetric tensor $T$ on $\RR ^n$ such that
$$
\eq{
X_i&(p)X_A=\indic_{(p,\cd)\not\in A}X_i(p)X_A+\sum_{j=1}^n \indic_{(p,j)\in
A}X_i(p)X_A\cr
&=\indic_{(p,\cd)\not\in A}X_{A\cup \{(p,i)\}}+\sum_{j=1}^n \indic_{(p,j)\in
A}X_i(p)X_j(p)X_{A\setminus \{(p,j)\}}\cr
&=\indic_{(p,\cd)\not\in A}X_{A\cup \{(p,i)\}}+\sum_{j=1}^n \indic_{(p,j)\in
A}\left(\d_{ij}+\sum_lT^{ij}_lX_l(p)\right)X_{A\setminus \{(p,j)\}}\cr
&=\indic_{(p,\cd)\not\in A}X_{A\cup \{(p,i)\}}+\indic_{(p,i)\in
A}X_{A\setminus (p,i)}+\sum_{j=1}^n \sum_{l=1}^n \indic_{(p,j)\in
A}T^{ij}_lX_{A\setminus \{(p,j)\}\cup \{(p,i)\}}\cr
}
$$
and we immediately recognize the formula for
$$
a^0_i(p)X_A+a^i_0(p)X_A+\sum_{p,l}T^{ij}_la^j_l(p)X_A.
$$\qed
\bigskip
Let us now return to quantum probabilistic structures and describe the
Fock space structure and its approximation by the atom chain. 

\spa{IV.}{Approximation of the Fock space by atom chains}

\def\pcc{\rP}
\def\rb{\RR}
\def\cb{\CC}
\def\pcn{\rP_n}
\def\fc{\rF}
\def\vi{\emptyset}

We recall the structure of the bosonic Fock space $\F$ and its basic
structure (see {\it e.g.} [At3] or [Pau] for details).

Let $\F=\G_s(L^2(\Rp;\CC^n))$ be the symmetric (bosonic) Fock space
over the space $L^2(\Rp;\CC^n)$.
We shall here give a very efficient presentation of that space, the
so-called {\it Guichardet interpretation} of the Fock space. 

Let $I=\{1,\ld,n\}$ and $\pcc$ be the set of finite subsets
$\{(s_1,i_1),\ld,(s_k,i_k)\}$ of
$\rb^+\times I$ such that the $s_i$ are mutually distinct. Then $\pcc =\cup_k \rP(k)$ where $\rP(k)$ is the set of $k$-elements subsets of
$\ \rb^+\times I$. By ordering the $\Rp$-part of the elements of
$\s\in\rP(k)$, the set $\rP(k)$ can be identified to the increasing simplex
$\Sigma _k = \{0<t_1 < \cdots <t_k\}\times I$ of $\rb^k\times I$. Thus $\rP(k)$
inherits a measured space structure from the product of Lebesgue measure on
$\ \rb^k$ and the counting measure on $I$. This also gives a measure structure on $\pcc$ if we specify
that on $\pcc_0 = \{\vi\}$ we put the measure $\delta
_\vi$. Elements of $\pcc$ are often denoted by $\sigma $, the measure on
$\pcc$ is denoted by $d\sigma $. The $\sigma $-field obtained this way on
$\pcc$ is denoted $\fc$.

We identify any element $\s\in\rP$ with a family
$\{\s_1,\ld,\s_{n}\}$ of (two by two disjoint)  subsets of $\Rp$ where 
$$
\s_i=\{s\in\Rp; (s,i)\in\s\}.
$$
For a $s\in\Rp$ we denote by $\{s\}_i$ the element
$\s=\{\emptyset,\ld,\emptyset,\{s\},\emptyset,\ld\emptyset\}$ of
$\rP$ (where $\{s\}$ is at the $i$-th position.

\def\Unn{\indic}

The {\it Fock space\/} $\Phi$ is the space $L^2(\pcc,\fc,d\sigma
)$. An element $f$ of $\Phi$ is thus a measurable function $f:\pcc \to
\cb$ such that
$$
\normca{f} = \int_{\pcc} |f(\sigma )|^2\ d\sigma  < \infty~.
$$
One can define, in the same way, $\pcc_{[a,b]}$ and $\Phi_{[a,b]}$ by
replacing $\rb^+$ with $[a,b]\subset \rb^+$. There is a natural
isomorphism between $\Phi_{[0,t]} \otimes  \Phi_{[t,+\infty[}$ given by
$h\otimes g \mapsto f$ where $f(\sigma ) = h(\sigma \cap [0,t]) \, g
(\sigma  \cap (t,+\infty[)$. 
Define also $\indic$ to be the {\it vacuum vector}, that is,
$\indic(\s)=\d_\emptyset(\s)$.

Define $\chi^i _t {\in} \Phi$ by
$$
\chi _t(\sigma ) = \cases {
\Unn_{[0,t]}(s) &if~~$\sigma =\{s\}_i$\cr 0 &otherwise. }
$$
\def\nb{\NN}
Then $\chi_t$ belongs to $\ \Phi_{[0,t]}$. We even have $\chi^i_t -
\chi^i_s \,\in \, \Phi_{[s,t]}$ for all $s\leq t$. This last property allows
to 
define a so-called {\it It\^o integral\/} on $\ \Phi $. Indeed, let $(g^i_t)_{t\ge
0}$ be  families in $\ \Phi $, for $i=1,\ld,n$, such that
\smallskip
 i) $t\mapsto \|g^i_t\|$ is measurable,
\smallskip
ii) $g^i_t \,{\in}\, \Phi _{[0,t]}$ for all $t$,
\smallskip
iii) $\int^\infty _0 \|g^i_t\|^2\ dt < \infty $,
\smallskip
\noindent then one defines $\sum_i\int^\infty _0 g^i_t\ d\chi^i_t$ to be the
limit in $\ \Phi $ of
$$
\sum_{i=1}^n \sum_{j=0}^\infty {1\over t_{j+1}-t_j} \int^{t_{j+1}}_{t_j} P_{t_j} g^i_s\ ds
\otimes \left(\chi^i _{t_{j+1}} - \chi^i_{t_j}\right)\eqno(3)
$$
where $P_t$ is the orthogonal projection onto $\Phi _{[0,t]}$ and
$\{t_j,~j{\in} \nb\}$ is a partition of $\ \rb^+$ which is understood to
be refining and to have its diameter tending to $0$. Note that
${1\over t_{j+1}-t_j} \int^{t_{j+1}}_{t_j} P_{t_j} g_s\ ds$ belongs to
$\ \Phi _{[0,t_j]}$, which explains the tensor product symbol in (3).

We get that $\sum_i\int^\infty _0 g^i_t \ d\chi^i_t$ is an element of $\Phi $
with
$$
\Big\|{\sum_i\int^\infty _0 g_t \ d\chi _t}\Big\|^2 = \sum_i\int^\infty _0
\ab{g^i_t}^2\,dt~.\eqno(4)
$$
Let $f {\in} L^2(\pcn)$; one can easily define the {\it iterated It\^o
integral\/} on $\ \Phi $.
$$
I_n(f) = \int_{\rP_n}f(\s)\, d\chi^{i_1}_{t_1}\ld d\chi^{i_n}_{t_n}
$$
by iterating the definition of the It\^o integral. We use the following notation: 
$$
I_n(f) = \int_{\pcn} f(\sigma )\ d\chi _\sigma
$$
which we extend, in an obvious way, to any $f\in\F$. 
We then have the following important representation.

\th{6.}{\it Any element $f$ of $\Phi $ admits an {\rm
abstract chaotic representation}
$$f = \int_{\pcc} f(\sigma )\ d\chi _\sigma $$
with
$$\|f\|^2 = \int_{\pcc} |f(\sigma )|^2\ d\sigma $$
and an {\rm abstract predictable representation}
$$f = f(\vi)\Unn + \sum_i\int^\infty _0 D^i_t f\ d\chi^i _t$$
with
$$\|f\|^2 = |f(\vi)|^2 + \sum_i\int^\infty _0 \|D^i_sf\|^2\ ds$$
where $[D^i_sf] (\sigma ) = f(\sigma  \cup \{s\}_i) \Unn_{\sigma \subset [0,s[}$.
}
\bigskip

Let us now recall the definitions of the basic noise operators
$a^i_j(t)$, $i,j=0,\ld, n$, on $\Phi $. They
are respectively defined by
$$
\eqalignno{
[a^0_i(t)f](\sigma ) &= \sum_{s{\in} \sigma_i \cap [0,t]} f(\sigma \setminus \{s\}_i),\cr
[a^i_0f](\sigma ) &= \int^t_0 f(\sigma \cup  \{s\}_i)\ ds,\cr
[a^i_jf](\sigma ) &= \sum_{s{\in} \sigma_i \cap [0,t]}\
f(\sigma\setminus\{s\}_i\cup\{s\}_j )\cr
}
$$
for $i,j\not=0$ and $$a^0_0(t)=tI.$$
\bigskip
There is a good common domain to all  these operators, namely
\def\dc{\rD}
$$\dc = \Big\{ f{\in} \Phi ~;~~\int_{\pcc} |\sigma | ~|f(\sigma )|^2\
d\sigma  < \infty  \Big\}~.$$

\def\scc{\rS}
Let $\ \scc = \{0 = t_0 < t_1 < \cdots <t_p < \cdots \}$ be a partition of
$\ \rb^+$ and $\delta (\scc) = \sup_i |t_{i+1}-t_i|$ be the diameter
of $\ \scc$. For fixed $\ \scc$, define $\ \Phi _p = \Phi
_{[t_p,t_{p+1}]}$, $i{\in} \nb$. We then have $\ \Phi  \simeq \otimes_{p{\in}
\nb} \Phi _p$ (with respect to the stabilizing sequence $(\Unn)_{p{\in}
\nb}$).

For all $p{\in} \nb$, define for $i,j\not=0$
$$
\eqalign{
X^i(p) &= {\chi^i _{t_{p+1}}-\chi^i _{t_p}\over \sqrt {t_{p+1}-t_p}} \in \Phi
_p~,\cr
a^i_0(p) &= {a^i_0(t_{p+1})-a^i_0({t_p})\over \sqrt{t_{p+1}-t_p}} P_{1]}~,\cr
a^i_j(p) &= P_{1]}\big(a^i_j(t_{p+1})-a^i_j(t_p) \big)P_{1]}~,\cr
a^0_j(p) &= P_{1]} {a^0_j(t_{p+1})-a^0_j({t_p})\over
\sqrt{t_{p+1}-t_p}}~,\cr}$$
where $P_{1]}$ is the orthogonal projection onto $L^2(\pcc_{ 1})$ and
where the above definition of $a^0_i(p)$ is understood to be 
valid on $\ \Phi _p$ only, with $a^0_i(p)$ being the identity operator $I$
on the others $\ \Phi _q$'s (the same is automatically true for $a^i_0$,
$a^i_j$). We put $a^0_0(p)=I$. 

\prp{7.}{\it We have
$$\eqalign{
&\cases{a^i_0(p)X^j(p)=\d_{ij}\Unn\cr
a^i_0 \Unn =0}\cr
\noalign{\vskip3pt}
&\cases{a^i_j(p)X^k(p)=\d_{ik}X^j(p)\cr
a^i_j \Unn =0}\cr
\noalign{\vskip3pt}
&\cases{a^0_j(p)X^i(p)=0\cr
a^0_j(p) \Unn =X^j(p).}\cr}
$$
}\qed
\bigskip
Thus the action of the operators $a^i_j$ on the $X^i(p)$ is
similar to the action of the corresponding operators on the atom chain
of section two. We are now going to construct the atom chain  inside $\Phi$.

We are still given a fixed partition $\scc$. Define $\TF (\scc)$ to
be the space of vectors $f {\in} \Phi $ which are of the form
$$
f =  \sum_{A{\in} \rP_N} f(A) X_A
$$
(with $\|f\|^2 = \sum_{A{\in} \rP_N} |f(A)|^2 < \infty $).

The space $\TF (\scc)$ is thus clearly identifiable to the atom chain
$\TF$; the operators $a^i_j(p)$ act on $\TF (\scc)$ exactly in the same way as the
corresponding operators on $\TF$. We have completely embedded the toy
Fock space into the Fock space.

Let $\ \scc = \{0=t_0 < t_1 <\cdots <t_p <\cdots\}$ be a fixed partition
of $\ \rb^+$. The space $\TF (\scc)$ is a closed subspace of $\ \Phi
$.  We denote by
$P_\rS$ the operator of orthogonal
projection from $\ \Phi $ onto $\TF (\scc)$. 
\bigskip
We are now going to prove that the Fock space $\Phi$ and its basic
operators $a^i_j(t)$ can be approached by the toy Fock
spaces $\TF(\scc)$ and their basic operators $a^i_j(p)$.

We are given a sequence $(\rS_p)_{p{\in} \NN}$ of partitions which are
getting finer and finer and whose diameter $\delta (\rS_p)$ tends to
$0$ when $p$ tends to $+\infty $. Let  $\TF(p) = \TF(\rS_p)$ and let $P_p$ be the orthogonal projector
onto $\TF(\rS_p)$,  for all $p{\in} \nb$.

\th{8.}{\it
\smallskip
i) For every $f {\in} \Phi$ there exists a sequence
$(f_p)_{p{\in} \NN}$ such that $f_p {\in} \TF(p)$, for all $p{\in} \NN$,
and $(f_p)_{p{\in} \NN}$ converges to $f$ in $\ \Phi$. 
\smallskip
ii) For all $i,j$
let 
$$
\e_{ij}=\frac12(\d_{0i}+\d_{0j}).
$$
 If $\scc_p = \{0 = t^p_0 < t^p_1 < \cdots < t^p_k <
\cdots \}$, then for all $t {\in} \rb^+$, the operators
$$
\sum_{k;t^p_k\leq t} {(t^p_{k+1}-t^p_k)}^{\e_{ij}} a^i_j(k)
$$
converge strongly on $\dc$ to $a^i_j(t)$.
}

\prf

{\it i)}~~As the $\ \scc_p$ are refining then the  $(P_p)_p$ form
an increasing family of orthogonal projection in $\Phi$. 
 Let $P_\infty  = \vee_p P_p$. Clearly,
for all $s\leq t$, all $i$ we have that $\chi^i_t - \chi^i_s$ belongs
 to Ran$P_\infty$. But by the construction of the It\^o integral 
and by Theorem 5, we have that the $\chi^i_t-\chi^i_s$ generate
$\Phi$. Thus $P_\infty  = I$. Consequently if $f {\in} \Phi$, the
sequence $f_p = P_p f$ satisfies the statements.

\medskip
{\it ii)} The convergence of 
$\sum_{k,t^p_k\leq t} {(t^p_{k+1} - t^p_k)}^{\e_{ij}} a^i_j(k)$ to
$a^i_j(t)$ is clear from the definitions when $i\not =0$. Let us check the
case of $a^0_i$. We have, for $f {\in} \dc$
$$
\bigg[\sum_{k;t^p_k\leq t} \sqrt{t^p_{k+1} - t^p_k} a^0_i(k) f\bigg]
(\sigma ) = \sum_{k;t^p_k\leq t} \Unn_{|\sigma \cap [t^p_k,t^p_{k+1}]|=1}
\sum_{s{\in} \sigma \cap [t^p_k,t^p_{k+1}]} f(\sigma \setminus \{s\}).
$$
Put $t^p = \inf \big\{ t^p_k{\in} \scc_p~; t^p_k \ge t\big\}$. We have
$$
\eqalign{
\Big\|& \sum_{k;t^p_k\leq t} \sqrt{t^p_{k+1} - t^p_k} a^0_i(k)  - a^0_i(t) f
\Big\|^2\cr
&= \int_\pcc \Big|\sum_{k;t^p_k\leq t} \Unn_{|\sigma \cap
[t^p_k,t^p_{k+1}]|=1} \sum_{s{\in} \sigma \cap [t^p_k,t^p_{k+1}]}
f(\sigma \setminus \{s\}) - \sum_{s{\in} \sigma \cap [0,t]} f(\sigma \setminus \{s\})
\Big|^2\ d\sigma \cr
&\leq 2 \int_{\pcc} \Big|\sum_{s{\in} \sigma \cap [t,t^p]} f(\sigma \setminus \{s\})
\Big|^2\ d\sigma \cr &\qquad +  2 \int_{\pcc} \Big|\sum_{k;t^p_k\leq t}
\Unn_{|\sigma \cap[t^p_k,t^p_{k+1}]|\ge 2} 
 \times \sum_{s{\in} \sigma \cap [t^p_k,t^p_{k+1}]} f(\sigma \setminus \{s\})
\Big|^2\ d\sigma.}
$$
\noindent For any fixed $\sigma $, the terms inside each of the integrals above
converge to $0$ when $p$ tends to $+\infty $. Furthermore we have, for
large enough $p$ ,
$$
\eqalign{
\int_{\pcc} \Big|\sum_{s{\in} \sigma \cap [t,t^p]} f(\sigma \setminus
\{s\})\Big|^2\ d\sigma
&\leq \int_{\pcc} |\sigma |\sum_{s{\in} \sigma\atop s\leq t+1} |f(\sigma \setminus
\{s\})|^2\ d\sigma \cr
&= \int^{t+1}_0 \int_{\pcc} (|\sigma |+1) |f(\sigma )|^2\ d\sigma \
ds\cr
&
\leq (t+1) \int_{\pcc} (|\sigma |+1) |f(\sigma )|^2\ d\sigma \cr}
$$
 which is finite for $f {\in} \dc$;
$$
\eqalign{
\int_{\pcc} \Big|\sum_{k;t^p_k\leq t}& \Unn_{|\sigma
\cap[t^p_k,t^p_{k+1}]|\ge 2}  \sum_{s{\in} \sigma \cap [t^p_k,t^p_{k+1}]}f(\sigma \setminus
\{s\})\Big|^2\ d\sigma\cr
&\leq \int_{\pcc} \Big(\sum_{k;t^p_k\leq t} \Unn_{|\sigma
\cap[t^p_k,t^p_{k+1}]|\ge 2}\  \Big| \sum_{s{\in} \sigma \cap [t^p_k,t^p_{k+1}]}f(\sigma \setminus
\{s\})\Big|\Big)^2\ d\sigma\cr
&\leq \int_{\pcc} \Big(\sum_{k;t^p_k\leq t}\ \sum_{s{\in} \sigma \cap
[t^p_k,t^p_{k+1}]} |f(\sigma \setminus
\{s\})|\Big)^2\ d\sigma\cr
&= \int_{\pcc} \Big(\sum_{s{\in} \sigma \atop s\leq t^p}  |f(\sigma \setminus
\{s\})|\Big)^2\ d\sigma\cr
&= \int_{\pcc}|\sigma |  \sum_{s{\in} \sigma \atop s\leq t^p}  |f(\sigma \setminus
\{s\})|^2\ d\sigma\cr
&\leq (t+1) \int_{\pcc} (|\sigma |+1) \big|f(\sigma )\big|^2\ d\sigma
\cr
}
$$
in the same way as above. So we can apply Lebesgue's theorem. This
proves {\it ii)}.

\findem

\spa{V.}{Multidimensional structure equations}
Let us recall some basic facts about normal martingales in $\RR^n$;
except for Theorem 13, all the statements in this section are taken from [A-E].

In the same way as the Fock space $\Phi  = \Gamma (L^2(\Rp;\CC))$ admits
probabilistic interpretations in terms of one-dimensional normal
martingales (see [At3]), the multiple Fock space $\Phi  = \Gamma (L^2(\Rp;\cb^n))$
admits probabilistic interpretations in terms of multidimensional
normal martingales. The point here is that the extension of the notion
of normal martingale, structure equation\dots to the
multidimensional case is not so immediate. Some interesting algebraic
structures appear.

A martingale $X= (X^1, \ld, X^n)$ with values in $\RR^n$
is called {\it normal} if $X_0 = 0$ and if, for all $i$ and $j$, the
  process $X^i_t X^j_t - \delta _{ij} t$ is a martingale. This is
  equivalent to saying that
$$
\langle X^i,X^j\rangle_t = \delta _{ij} t
$$
for all $t \in \Rp$, or else this is equivalent to saying that the
process
$$
[X^i,X^j]_t - \delta _{ij} t
$$
is a martingale.

A normal martingale $X= (X^1 \ld X^n)$ in $\RR^n$ is
said to {\it satisfy a structure equation} if each of the martingales
$[X^i,X^j]_t - \delta _{ij}t$ is a stochastic integral with respect to
$X$:
$$
[X^i,X^j]_t = \delta _{ij}t + \sum^n_{k=1} \int^t_0 T ^{ij}_k
(s)\ dX^k_s
$$
where the $T ^{ij}_k$ are predictable processes.

Any family $\left\{A^{ij}_k~;~i,j,k \in \{1\ld n\}\right\}$  of real
numbers is identified to a 3-tensor, that is, a linear map $A$ from $\RR^n$ to $\RR^n \otimes \RR^n$
by
$$
(Ax)_{ij} = \sum^n_{k=1} A^{ij}_k x_k~.
$$
Such a family is said to be {\it diagonalizable in some orthonormal
  basis} if there exists an orthonormal basis $\{e^1\ld e^n\}$ of
$\RR^n$ for which
$$
Ae^k = \lambda _k e^k \otimes e^k
$$
for all $k = 1 \ld n$ and for some {\it eigenvalues} $\lambda _1\ld
\lambda _n \in \RR$.

A family $\left\{A^{ij}_k~;~i,j,k\in \{1\ld n\}\right\}$ is called
{\it doubly symmetric} if
\smallskip
i) $(i,j,k) \mapsto A^{ij}_k$ is symmetric on $\{1 \ld
  n\}^3$ and
\smallskip
ii) $(i,j,i',j') \mapsto \sum^n_{k=1} A^{ij}_k
  A^{i'j'}_k$ is symmetric on $\{1\ld n\}^4$.

\th{9.}{\it
For a family $\left\{A^{ij}_k~;~i,j,k\in \{1\ld n\}\right\}$ of real
numbers, the following assertions are equivalent.
\smallskip
 i) $A$ is doubly symmetric.
\smallskip
 ii) $A$ is diagonalizable in some orthonormal basis.
}
\bigskip
This means that the condition of being doubly symmetric is the exact
extension to 3-tensors of the symmetry property for matrices
(2-tensors): it is the necessary and sufficient condition for being
diagonalisable in some orthonormal basis.
\bigskip
A family $\{x^1 \ld x^k\}$ of elements of $\RR$ is called
{\it orthogonal family} if the $x^i$ are all different from $0$ and
are two by two orthogonal.

\th{10.}{\it
There is a bijection between the doubly symmetric families $A$ of
$\RR^n$ and the orthogonal families $\Sigma $ which is given by
$$
Af = \sum_{x\in \Sigma } \frac{1}{\|x\|^2} \langle x,f\rangle
\, x\otimes x
$$
and
$$
  \Sigma  = \left\{ x\in \RR^n\setminus\{0\}~;~ Ax = x\otimes x\right\}.
$$}

\bigskip

These algebraic preliminaries are used to determine the
behaviour of the multidimensional normal martingales.

\th{11.}{\it
Let $X$ be a normal martingale in $\RR^n$ satisfying a structure
equation
$$
[X^i,X^j]_t = \delta _{ij}t + \sum^n _{k=1} \int^t_0 T
^{ij}_k (s)\ dX^k_s.
$$
Then for \aa $(t,\omega )$ the family $\{T ^{ij}_k(s,\omega
)~;~i,j,k = 1\ld n\}$ is doubly symmetric. If $\Sigma _t (\omega )$ is
its associated  orthogonal system and if $\pi _t(\omega )$ denotes the
orthogonal projection onto $(\Sigma _t(\omega ))^\perp$, then the
continuous part of $X$ is given by
$$
X^{c,i}_t = \sum^n_{j=1} \int^t_0 \pi ^{ij}_s\ dX^j_s;
$$
the jumps of $X$ happen only at totally inaccessible times and
they satisfy
$$
\Delta  X_t(\omega ) \in \Sigma _t(\omega )~.
$$
}

We can now study a basic example. The simplest case occurs when
$T $ is constant in $t$. Contrarily to the unidimensional case,
this situation is already rather rich.

\prp{12.}{\it
Let $T $ be a doubly symmetric family on $\RR^n$. Let $\Sigma $ be
its associated orthogonal system. Let $ B$ be a Brownian motion
with values in the Euclidian space $\Sigma ^\perp$. For each $x\in
\Sigma $, let $N^x$ be a Poisson process with intensity
$\|x\|^{-2}$. We assume $ B$ and all the $N^x$ to be
independent. Then the martingale
$$
X_t = B_t + \sum_{x \in \Sigma } (N^x_t - \|x\|^{-2}_t) x
$$
satisfies the constant coefficient structure equation
$$
  [X^i,X^j]_t = \delta _{ij}t + \sum^n_{k=1} T^{ij}_k  \ X^k_t~.
$$
Conversely, every normal martingale which is solution of the above
equation has the same law as $ X$.
}
\bigskip
Finally, let us recall a particular case of a theorem proved in [At2],
which has the advantage of not needing the introduction of quantum
stochastic integrals and of being sufficient for our purpose.

\th{13.}{\it Let $X$ be a normal martingale in $\RR^n$ which satisfies
a structure equation of the above form :
$$
  [X^i,X^j]_t = \delta _{ij}t + \sum^n_{k=1} T^{ij}_k  \ X^k_t~.
$$
Then $(X_t)_t$ possesses the chaotic representation property. Furthermore, the  space
$L^2(\Omega, \rF, P)$, where $(\Omega, \rF, P)$ is the canonical space
associated with $(X_t)_t$, is naturally isomorphic
to  $\F$, by identification of the chaotic expansion of $f$ with the
element $\wt f$ of $\F$ whose abstract chaotic expansion has the same
coefficients. 

Within this identification the operator of
multiplication by $X^k_t$ is equal to
$$
\rM_{X^k_t}=a^0_k(t)+a^k_0(t)+  \sum^n_{i,j=1} T^{ij}_k  a^i_j(t).
$$}

\spa{VI}{Convergence to normal martingales}
Now we can close the circle under the form of a kind of commutative
diagram and establish some convergence theorem.

Starting from an obtuse random random variable $X$ depending on a
parameter $h\in\Rp$, with associated
sesqui-symmetric tensor $T$, we associate a random walk
$(X_p)_{p\in\NN}$of i.i.d. random variables with the same law as
$X$. By Theorem 2, the renormalized random walk 
$$\wt X_n=\sqrt h\,X_n$$
satisfies 
the discrete time structure equation
$$
\wt X\otimes\wt X=hI+\wt T(\wt X)
$$
where $\wt T^{ij}_k=\sqrt h\,T^{ij}_k$. The tensor $\wt T$ is also
sesqui-symmetric but with the relation
\smallskip
ii')\ \ $(i,j,l,m)\longmapsto\sum_kT^{ij}_kT^{lm}_k+h\d_{ij}\d_{lm}$ is
symmetric.
\smallskip
Theorem 5 shows that the associated multiplication operator by $\wt X$
is given by  
$$
U\rM_{\wt X_i(k)}U^\ast=\sqrt h (a^0_i(k)+a^i_0(k))+\sum_{j,l=1}^n \wt
T^{jl}_i\,a^j_l(k).
$$
By Proposition 7 we can embed this situation inside the Fock space
$\F$ and we
get a family of operators on $\F$ such that 
$$
\sum_{k\leq [t/h]} U\rM_{\wt X_i(k)}U^\ast
$$
converges strongly on $\rD$ to
$$
X_t=a^0_i(t)+a^i_0(t)+\sum_{j,l=1}^n S^{jl}_i a^j_l(t)
$$
where $S^{jl}_i=\lim_{h\rightarrow 0} \wt T^{jl}_i$, by Theorem
8. Because of the relation ii') above, the limit tensor $S$ is
automatically doubly-symmetric.

Thus by Theorem 13, the operators $X_t$ are the canonical
multiplication operators by a normal martingale, solution of the
structure equation 
$$
  [X^i,X^j]_t = \delta _{ij}t + \sum^n_{k=1} S^{ij}_k  \ X^k_t~.
$$
\bigskip
From the above we  see that only  the coefficients $\wt T^{ij}_k$
which admit a limit $S^{ij}_k$, when $h\rightarrow 0$, contribute to
the limit normal martingale $\pro X$. This means that only the coefficients
$T^{ij}_k$ which have a dominant term of order $1/\sqrt h$ will
contribute non-trivialy to the limit. A smaller dominant term gives 0
in the limit and a larger dominant term will not admit a limit.

If the obtuse random variable $X$ is given one direction for which its
probability is of order $h$, then, by Proposition 1 iv), the length of
the jump in that direction is of order $1/\sqrt h$.  The associated
tensor will then get terms $T^{ij}_k$ of order $1/\sqrt h$ too
(Theorem 3). Thus in the limit this terms will participate to the
tensor $S$. By Proposition 12, these terms $S^{ij}_k$ will participate
to the Poisson-type behaviour of the normal martingale.

In the same way one gets easily conviced that the directions of $X$
which are visited with a probability of constant order, or of bigger
order than $h$ will contribute to the diffusive part of the
martingale.
\bigskip
Note that, in order to understand the above discussion in
probabilistic terms it is not necessary to pass throught the
representation in terms of creation and annihilation operators. One
can directly approach a normal martingale in $\RR^n$ by some obtuse
random walks (this has be achieved explicitely in [Tav]). But this was
not our purpose here to detail this approximation. We just wanted to
show how it is naturally related to the approximation of the Fock
space by state spaces of $(n+1)\,$-level atom chains.
\bigskip
We have already a convergence of the random walk to a
normal martingale of which the law is given by Theorem 12. Yet this strong convergence of multiplication
operators is not easy to translate into probabilistic language,
because determining which random variables in $L^2(\Omega, \rF, P)$
are sent to $\rD$ by identification is not an easy problem (it amounts
to studying the integrability properties of the chaotic expansion of
random variables).

Obtaining the convergence in law in the above framework requires some
more work. This is what the next theorem does.

\th{14.}{\it With the above notations, the random variable 
$$ \sqrt h \, \sum_{n=1}^{[t/h]} X_n$$
converges to $X_t$ in law, for almost all $t$.}

\prf
\def\a{\alpha}
Developping all the details of this proof would need much more serious and
 longer developments which are not compatible with the length and the
spirit of this article. This is the reason why we adopt a more concise
style in the following proof.
\bigskip
Choose $\a$ in $\RR^n$ with coordinates $\a_1, \ld, \a_n$ and denote
by $M(\a,h,k)$ the 
operator of multiplication by 
$$ \langle \a , \sqrt h \, \sum_{p=1}^k X_p\rangle_{\RR^n}$$
on $L^2(A^{\otimes \NN},\rA^{\otimes \NN}, P^{\otimes\NN})$. Then denote by $u_k$ the operator $\exp(i\,M(\a,h,k))$. The family $u_k$
satisfies the equations 
$$ \cases{u_{k+1}= e^{i\langle\a, \sqrt h X_{k+1}\rangle} u_k\cr
\ u_0= I.}$$
Now, by a slight adaptation of Theorem 19 in [A-P] (the terms $D_{ij}$
have to be allowed to converge as $h$ goes to zero instead of being
constant, but this is immediate from the proof), the sequence
$u_{[t/h]}$, seen as a sequence of operators on $\Phi$, converges
strongly to $U_t$ for almost all $t$, where $(U_s)$ solves
$$\cases{dU_s = \sum_{j,l = 0}^n L^j_l U_s \, da^j_l(s) \cr \ U_0 =I}$$
with 
$$\eqalign{L^0_0 &= -\frac12 W^\ast W \cr L^0_l &= W_l \cr L^j_0 &-(W^\ast S)_j \cr L^j_l &= S^j_l - I}$$
where
$$\eqalign{S &= \exp(iD) \cr W&=(S-I)/D \, \a}$$
and $D$ is the $n\times n$ matrix with coefficients
$$ D_{jl}= \sum_{i=1}^n \a_i T^{jl}_i.$$
Note that we do not assume that $D$ is invertible; if it is not then
$(S-I)/D$ is still defined by its series expansion.

Such an equation, called a Hudson-Parthasarathy equation, always has a
solution made of unitary operators $U_s$, and such an equation is
unique (for such results we refer the reader to [Par]).
\smallskip

On the other hand, consider the operator $V_t$ which is the value for
$s=1$ of the unitary semigroup $e^{isM(\a,t)}$ where $M(\a,t)$ is the
operator of multiplication by $X_t$ (on the space $L^2(\Omega)$
identified with $\Phi$). Remember that $M(\a,t)$ has the
representation 
$$ \sum_{i=1}^n \a_i \Big( a^0_i(t) + a^i_0(t) + \sum_{j,l=1}^n
T^{jl}_i \, a^j_l(t) \Big)$$
on $\rD$.
We now wish to apply Vincent-Smith's formula (see [ViS])
for $e^{i M(\a,t)}$, but Vincent-Smith's result would require $M(\a,t)$ to
belong to the class of regular semimartingales (as in [At4]);
here it would need the additional property of being bounded.

Yet in our case it is easy to prove {\it a posteriori} that the
formula holds; this is done by an usual trick. Checking
that a formal application of Vincent-Smith's formula for $\exp(is M(\a,t))$
gives a quantum stochastic integral which is a bounded operator,
a strongly continuous semigroup of the parameter $s$, and that it is
unitary, its generator can then be computed  and shown to be equal
$M(\a,t)$ on some good 
domain ({\it e.g.} the space of coherent vectors). Since $M(\a,t)$, as
a linear combination of the fundamental 
operators $a^j_l(t)$, is known to be essentially selfadjoint on the
exponential domain, Stone's theorem proves the validity of the
integral representation on $\Phi$ for all $\exp(is M(\a,t))$.

Let us therefore apply Vincent-Smith's formula for $\exp(is
M(\a,t))$. It simplifies greatly since the operator $\langle \a,
X_t\rangle_{\RR^n}$, as a scalar multiplication, commutes with all
other coefficients. We obtain therefore the equation
$$ dV_s = \sum_{j,l=1}^n K^j_l V_s\, da^j_l(s)$$
with 
$$ \eqalign{
K^0_l &= \int_0^1 \exp\big(i(1-u)D\big) \a \, du \cr 
K^j_0 &= \int_0^1 \a^\ast \exp\big(iuD\big) \, du \cr
K^j_l &= \int_0^1 \exp\big(i(1-u)D\big) D\, du \cr
K^0_0 &= \int_0^1\int_0^1 u \, \a^\ast \exp\big( i u (1-v) D\big) \a \,
du\, dv \cr
&= \int_0^1 \a^\ast \exp\big(iuD\big) \a\, du
}$$
and it is clear that the $K^j_l$ equal the corresponding
$L^j_l$. Since $V_0= U_0$ we obtain the equality $U_t=V_t$ for all
$t$.

\smallskip
On the other hand, for all $\a$, the operator $M(\a,h,k)$ is
selfadjoint on $\TF(h)$. By Nelson's theorem, there exists a dense set
of vectors $f$ such that 
$$ u_k f = \sum_{j\in\NN} \frac{i^j M(\a,h,k)^j}{j\,!} f.$$
When seen on the $L^2(\Omega, \rF, P)$, this means that almost
everywhere, the equality
$$ (u_k f)(\o) = \exp (i \langle \a,\sqrt h \, \sum_{p=1}^k X_p
(\o)\rangle_{\RR^n}) \, f(\o) $$
holds for all $f$ in the previously specified dense set. Note that this
equality is only apparently trivial because of our identification.

The almost everywhere equality can therefore be extended to all $f$ in
$\TF(h)$, all $k$. Similarly we can obtain the almost everywhere
equality
$$ (U_t f)(\o) = \exp(i \langle \a, X_t(\o)\rangle _{\RR^n})$$
for all $f$ in $\Phi$.

\smallskip
The above equalities and the strong convergence prove in particular
that
$$ \EE \big(\exp(i\langle \a, \sum_{p=1}^{[t/h]} X_p
\rangle_{\RR^n})\big)= \langle \O, u_{[t/h]} \O \rangle_{\TF(h)}$$
converges to 
$$ \langle \O, U_t \O \rangle_{\Phi} = \EE \big(\exp(i\langle \a, X_t
\rangle_{\RR^n}) \big)$$
as $h$ tends to zero. This holds for all $\a$, so that the conclusion
holds.

\findem

\spa{VII.}{Some approximations of $2$-dimensional noises}

We end this article by computing some simple and illustrative examples
in the case $n=2$. 
\bigskip
We consider, in the case $n=2$, an obtuse random variable $X$ which
takes the values $v_1=(a,0)$, $v_2=(b,c)$ and $v_3=(b,d)$ with
respective  probabilities $p,q,r$. In order that $X$ be obtuse we put
$$
a=\sqrt{1/p-1},\  b=-1/a,\  c=\sqrt{1/q-1-b^2},\  d=-\sqrt{1/r-1-b^2}.
$$
Let us call $S$ this set of values for $X$ and $p_s$ the probability
associated to $s\in S$. The associated sesqui-symmetric 3-tensor $T$ is given by
$$
T(v)=\sum_{s\in S}p_s \, <s,x>s\otimes s.
$$
\bigskip
For example, in the case $p=1/2$, $q=1/3$ and $r=1/6$ we get $a=1$,
$b=-1$, $c=1$ and $d=-2$. The tensor $T$ is then given by
$$
T(v)=\left(\matrix{0&-y\cr-y&-x-y\cr}\right)
$$
if $v=(x,y)$.
Thus the multiplication operator by $X_1$ is equal to
$$
X_1=a^1_0+a_1^0-a^2_2
$$
and the multiplication operator by $X_2$ is equal to
$$
X_2=a^2_0+a^0_2-(a^1_2+a^2_1+a^2_2).
$$
Now we consider a random walk $(X(k))_{k\geq0}$ made of independent
copies of this 
random variable $X$, with time step $h$. In the framework of the Fock
space approximation 
described above, the operator 
$$
\sum_{k;kh\leq t}\sqrt h\, X_1(k)
$$
converges, both in the sense of convergence of multiplication
operators and in law, to
$$
a^0_1(t)+a^1_0(t)
$$
and the operator  
 $$
\sum_{k;kh\leq t}\sqrt h\,X_2(k)
$$
converges to
$$
a^0_2(t)+a^2_0(t).
$$
This means that the limit process $X(t)$ is a 2-dimensional Brownian
motion. Indeed, the above representation shows that the associated
doubly-symmetric tensor $\Phi$ is null and thus $X$ satisfies the
structure equation 
$$
\eqalign{
d[X_1,X_1]_t&=dt\cr
d[X_1,X_2]_t&=0\cr
d[X_2,X_2]_t&=dt\cr
}
$$
which is exactly the structure equation verified by two independent
Brownian motions.
\bigskip
It is clear, that whatever the values of $p,q,r$ are, if they are
independent of the time step parameter $h$, we will always obtain a
2-dimensional Brownian motion as a limit of this random walk.
\bigskip
When some of the probabilities $p,q$ or $r$ depend on $h$ the
behaviour is very different. Let us follow two examples.
\smallskip
In the case $p=1/2$, $q=h$ and $r=1/2-h$ we get
$$
a=1,\ b=-1,\ c={1\over{\sqrt h}}+O(h^{1/2}),\ d=-2\sqrt h+o(h^{3/2}).
$$
For the tensor $T$ we get 
$$
T(v)=\left(\matrix{0+o(h^{5/2})&-y+o(h^2)\cr
-y+o(h^2)&-{y\over{\sqrt h}}-x+o(h^{1/2})\cr}\right).
$$
The multiplication operators are then given by
$$
X_1=a^1_0+a^0_1-a^2_2+O(h^2)
$$
and
$$
X_2=a^2_0+a^0_2-(a^1_2+a^2_1)+{1\over{\sqrt h}}a^2_2+O(h^{1/2}).
$$
In the same limit as above we thus obtain the operators
$$
a^1_0(t)+a^0_1(t)
$$
and
$$
a^2_0(t)+a^0_2(t)-a^2_2(t).
$$
This means that the coordinate $X_1(t)$ is a Brownian motion and
$X_2(t)$ is an independent Poisson process, with intensity 1 and
directed upwards.  Indeed, the associated tensor $\Phi$ is given by
$$
\Phi(v)=\left(\matrix{0&0\cr
0&-y}\right)
$$
and the associated structure equation is 
$$
\eqalign{
d[X_1,X_1]_t&=dt\cr
d[X_1,X_2]_t&=0\cr
d[X_2,X_2]_t&=dt+dX_2(t)\cr
}
$$
which is the structure equation of the process we described.
\bigskip
The last example we shall treat is the case $p=1-2h$, $q=r=h$. We get,
for the dominating terms
$$
a=\sqrt 2\sqrt h,\ b=-{1\over{\sqrt 2}}{1\over{\sqrt h}},\
c={{1}\over{\sqrt 2}}{1\over{\sqrt h}},\  d=-{1\over{\sqrt
2}}{1\over{\sqrt h}},
$$
and
$$
\eqalign{
X_1&=a^1_0+a^0_1-{1\over{\sqrt 2}}{1\over{\sqrt
h}}a^2_2+{{1}\over{\sqrt 2}}{1\over{\sqrt h}}a^1_1\cr
X_2&=a^2_0+a^0_2-{1\over{\sqrt 2}}{1\over{\sqrt
h}}(a^1_2+a^2_1).\cr
}
$$
The limit process is then solution of the structure equation
$$
\eqalign{
d[X_1,X_1]_t&=dt-{1\over{\sqrt 2}}\,dX_1(t)\cr
d[X_1,X_2]_t&=-{1\over{\sqrt 2}}\, dX_2(t)\cr
d[X_2,X_2]_t&=dt-{1\over{\sqrt 2}}\, dX_1(t).\cr
}
$$
The associated tensor is easy to diagonalise and one finds the
eigenvectors
$$
(-1/\sqrt 2,1/\sqrt 2)\ \ {\rm and}\ \ (-1/\sqrt 2,-1/\sqrt 2).
$$
The limit process is made of two independent Poisson processes, with
intensity 2 and respective direction (-1,1) and (-1,-1). 
\bigskip
\vfill \eject
\spa{}{References}
\bigskip\noindent\hangindent=1cm\hangafter=1
[At1]: Attal S.: ``{\it Approximating the Fock space with the toy
Fock space}'', { S\'eminaire de Probabilit\'es XXXVI} (2003),
Lect. Notes in Math. 1801, Springer verlag, p.\thsp 477-491
\bigskip\noindent\hangindent=1cm\hangafter=1
[At2]: Attal S.: ``{\it Semimartingales non commutatives et applications
aux endomorphismes browniens}'', Thesis of Strasbourg University, 1994.  
\bigskip\noindent\hangindent=1cm\hangafter=1
[At3]: Attal S.: ``{\it Classical and quantum stochastic calculus}'', {
Quantum Probability and Related Topics X} (1998), World Scientific,
p.\thsp 1-52. 
\bigskip\noindent\hangindent=1cm\hangafter=1
[At4]: Attal S.: ``{\it An algebra of non-commutative bounded
semimartingales --\allowbreak Squa\-re and angle quantum brackets}'', {Journal of
Functional Analysis} { 124} 
(1994), p.\thsp 292-332. 
\bigskip\noindent\hangindent=1cm\hangafter=1
[A-E]: Attal S. and  Emery M.:  ``{\it Equations de structure pour des
martingales vectorielles}'', { S\'eminaire de Probabilit\'es XXVIII}
(1994), Springer Verlag, p. 256-278. 
\bigskip\noindent\hangindent=1cm\hangafter=1
[A-P]: Attal S. and Pautrat Y.:  ``{\it From discrete to continuous
quantum interactions}'', {preprint submited to Duke Mathematical Journal}. 
\bigskip\noindent\hangindent=1cm\hangafter=1
[Gui]: Guichardet A.: ``{\it  
 Symmetric Hilbert spaces and related topics}'',  Lecture Notes in
Mathematics, Vol. 261.  
 Springer-Verlag, Berlin-New York, 1972.
\bigskip\noindent\hangindent=1cm\hangafter=1
[Par]: Parthasarathy K.R.: ``{\it An introduction to quantum stochastic
calculus}'', {\sl Mo\-no\-graphs in Mathematics 85}, Birkh\"auser (1992).
\bigskip\noindent\hangindent=1cm\hangafter=1
[Pau]: Pautrat Y.: ``{\it Des matrices de Pauli aux bruits quantiques}'',
{ Thesis of Grenoble University, 2003.}
\bigskip\noindent\hangindent=1cm\hangafter=1
[Tav]: Taviot G.: ``{\it Martingales et \'equations de structure :
\'etude g\'eom\'etrique}'',
{ Thesis of Strasbourg University, 1999.}
\bigskip\noindent\hangindent=1cm\hangafter=1
[ViS]:  Vincent-Smith G.F.: ``{\it The It\^o formula for quantum
semimartingales}'', {\sl Proc. London Math. Soc. (3)  75}  (1997),
no. 3, pp.671-720.
\bigskip\bigskip
\line{St\'ephane Attal\hf}
\line{ Institut Fourier, U.M.R. 5582\hf}

\line{ Universit\'e de Grenoble I, BP 74\hf }

\line{ 38402 St-Martin d'H\`eres cedex, France\hf }

\bigskip
\line{Yan Pautrat\hf}
\line{ McGill mathematics and statistics\hf}

\line{ 805 Sherbrooke West\hf }

\line{ Montreal, QC, H3A 2K6, Canada\hf }
\end